\documentclass{amsart}
\usepackage{amsmath,amssymb,amsfonts,amsthm}
\usepackage[dvips]{graphicx}

\begin{document}
\newcommand{\omegaone}{\ensuremath{\omega_1}}
\newcommand{\lomegaone}{\ensuremath{\mathcal{L}_{\omega_1,\omega}}}
\newcommand{\alephs}[1]{\ensuremath{\aleph_{#1}}}
\newcommand{\alephalpha}{\alephs{\alpha}}
\newcommand{\alephomegaone}{\alephs{\omegaone}}
\newcommand{\alephalphaplus}{\alephs{\alpha+1}}
\newcommand{\beths}[1]{\ensuremath{\beth_{#1}}}
\newcommand{\bethalpha}{\beths{\alpha}}
\newcommand{\bethomegaone}{\beths{\omegaone}}
\newcommand{\bethalphaplus}{\beths{\alpha+1}}
\newcommand{\z}{\ensuremath{\mathcal{Z}}}
\newcommand{\M}{\ensuremath{\mathcal{M}}}
\newcommand{\N}{\ensuremath{\mathcal{N}}}
\newcommand{\A}{\ensuremath{\mathcal{A}}}
\newcommand{\B}{\ensuremath{\mathcal{B}}}
\newcommand{\F}{\ensuremath{\mathcal{F}}}
\newcommand{\D}{\ensuremath{\mathcal{D}}}
\newcommand{\C}{\ensuremath{\mathcal{C}}}
\newcommand{\E}{\ensuremath{\mathcal{E}}}
\newcommand{\G}{\ensuremath{\mathcal{G}}}
\newcommand{\Umax}{\ensuremath{U_{max}}}
\newcommand{\tildeP}{\ensuremath{\tilde{P}}}
\newcommand{\mnkappafull}{$(\M,\N,\kappa)$-full\;}
\newcommand{\mneat}{\M-neat\;}
\newcommand{\mnneat}{(\M,\N)-neat\;}
\newcommand{\mnhappy}{(\M,\N)-happy\;}
\newcommand{\mnrich}{(\M,\N)-rich\;}
\newcommand{\mnfull}{(\M,\N)-full\;}
\newcommand{\mrich}{\M-rich\;}
\newcommand{\mfull}{\M-full\;}
\newcommand{\AminusMN}{\ensuremath{\A\setminus(\M\cup\N)\;}}
\newcommand{\lang}[1]{\ensuremath{\mathcal{L}_{#1}}}
\newcommand{\langhat}{\ensuremath{\widehat{\lang{}}}}
\newcommand{\phialpha}{\ensuremath{\phi_\alpha}}
\newcommand{\phikappa}{\ensuremath{\phi_\kappa}}
\newcommand{\phiM}{\ensuremath{\phi_{\M}}}
\newcommand{\phiF}{\ensuremath{\phi_{\F}}}
\newcommand{\philtok}{\ensuremath{\phi_{\ltok}}}
\newcommand{\phialphaplus}{\ensuremath{\phi_{\alpha+1}}}
\newcommand{\existslambda}{\ensuremath{\exists^{\omega\le\cdot\le\lambda}}}
\newcommand{\Malpha}{\ensuremath{\mathcal{M}_\alpha}}
\newcommand{\ch}{\ensuremath{\mathcal{CH}_{\omega_1,\omega}}}
\newcommand{\homch}{\ensuremath{\mathcal{HCH}_{\omega_1,\omega}}}
\newcommand{\ltok}{\ensuremath{\lambda^{\kappa}}}
\newcommand{\ktok}{\ensuremath{\kappa^{\kappa}}}
\newcommand{\ltoomega}{\ensuremath{\lambda^{\omega}}}
\newcommand{\kappaplus}{\ensuremath{\kappa^{+}}}
\newcommand{\kappaplusplus}{\ensuremath{\kappa^{++}}}
\newcommand{\lambdaalpha}{\ensuremath{\lambda_{\alpha}}}
\newcommand{\lambdaalphas}[1]{\ensuremath{\lambda_{\alpha_{#1}}}}
\newcommand{\dalpha}{\ensuremath{D_{\alpha}}}
\newcommand{\dalphaplus}{\ensuremath{D_{\alpha+1}}}
\newcommand{\dalphas}[1]{\ensuremath{D_{#1}}}
\newcommand{\taukappa}{\ensuremath{\tau_{\kappa}}}
\newcommand{\fuptox}[2]{\ensuremath{#1\upharpoonright_{(-\infty,#2)}}} 
\newcommand{\Sinfty}{\ensuremath{S_{\infty}}}
\newcommand{\Skappa}{\ensuremath{S_{\kappa}}}
\newcommand{\f}[1]{\ensuremath{f(\{#1\})}}
\newcommand{\cfun}[1]{\ensuremath{c(\{#1\})}}
\newcommand{\tk}{\ensuremath{t_k}}
\newcommand{\newlo}{\ensuremath{ \triangleleft}}
\newcommand{\simm}{\sim_m}
\newtheorem{theorem}{Theorem}[section]
\newtheorem{lemma}[theorem]{Lemma}
\newtheorem{corollary}[theorem]{Corollary}
\theoremstyle{df}
\newtheorem{df}[theorem]{df}
\newtheorem{claim}{Claim}
\newtheorem{subclaim}{Subclaim}
\newtheorem{note}{Note}
\theoremstyle{observation}
\newtheorem{observation}[theorem]{Observation}
\newtheorem{open}{Open Question}
\newtheorem{example}[theorem]{Example}
\newtheorem{xca}[theorem]{Exercise}
\theoremstyle{remark}
\newtheorem{remark}[theorem]{Remark}
\numberwithin{equation}{section}

\title{Linear Orderings and Powers of Characterizable Cardinals}

\author {Ioannis  Souldatos}
\address{Department of Mathematics,273 Wissink Hall, Minnesota State University, Mankato, MN, USA}
\email{ioannis.souldatos@mnsu.edu}



\date{November 3, 2010}
\begin{abstract}
The current paper answers an open question of  \cite{characterizablecardinals}.

We say that a countable model $\M$ characterizes an infinite cardinal $\kappa$, if the Scott sentence of $\M$ has a model in cardinality $\kappa$, but no models in cardinality $\kappa^+$. If $\M$ is linearly ordered by $<$, we will say that the linear ordering $(\M,<)$ characterizes $\kappa$, or that $\kappa$ is characterizable by $(\M,<)$.

From \cite{Knight} we can deduce that if $\kappa$ is characterizable, then $\kappa^+$ is characterizable by a linear ordering (see theorem \ref{loofsizekappa+}, corollary \ref{loforkappaplus}). From \cite{characterizablecardinals} we know that if $\kappa$ is characterizable by a dense linear ordering, then $2^\kappa$ is characterizable (see theorem \ref{loimplies2k}).

We show that if $\kappa$ is \emph{homogeneously} characterizable (cf. definition \ref{defhch}), then $\kappa$ is characterizable by a dense linear ordering, while the converse fails (theorem \ref{homchimpliesloch}).

The main theorems are: 1)  If $\kappa>2^{\lambda}$ is a characterizable cardinal, $\lambda$ is characterizable by a dense linear ordering and $\lambda$ is the least cardinal such that $\kappa^{\lambda}>\kappa$, then $\kappa^{\lambda}$ is also characterizable (theorem \ref{leastlambda}), 2) if $\alephalpha$ and $\kappa^{\alephs{\alpha}}$ are characterizable cardinals, then the same is true for $\kappa^{\alephs{\alpha+\beta}}$, for all countable $\beta$ (theorem \ref{powerclusters}).

Combining these two theorems we get that if $\kappa>2^{\alephs{\alpha}}$ is a characterizable cardinal, $\alephs{\alpha}$ is characterizable by a dense linear ordering and $\alephs{\alpha}$ is the least cardinal such that $\kappa^{\alephs{\alpha}}>\kappa$, then for all $\beta<\alpha+\omegaone$ $\kappa^{\alephs{\beta}}$ is characterizable (theorem \ref{powerclusters2}). Also if $\kappa$ is a characterizable cardinal, then $\kappa^{\alephs{\alpha}}$ is characterizable, for all countable $\alpha$ (corollary \ref{countable powers}).
\end{abstract}

\maketitle


\section{Structure of the paper}
Throughout the whole paper we work with countable languages $\lang{}$ and when we refer to a dense linear ordering we mean a dense linear ordering without endpoints. The first two sections provide some background material for the characterizable cardinals and for the dense linear orderings respectively. Section \ref{kappatoalephaone} contains the construction that proves the following

\begin{theorem} If $\kappa$ is a characterizable cardinal, then $\kappa^{\alephs{1}}$ is also a characterizable cardinal.
\end{theorem}
This appears as theorem \ref{kappatoalpheoneinch} in section \ref{kappatoalephaone} and it will be easily generalized to $\lambda\ge\alephs{1}$ in the last section.

\section{Characterizable cardinals}
 This section provides the necessary background on  characterizable cardinals.

\begin{df}\label{charcarddef} We say that a $\lomegaone$ sentence $\phi$ \emph{characterizes}
\alephalpha, or that $\alephalpha$ is\emph{ characterizable}, if $\phi$ has models in cardinality \alephalpha, but not in cardinality
\alephalphaplus. If $\phi$ is the Scott sentence of a
countable model $\M$ (or $\phi$ is any complete sentence), we say that $\M$
\emph{completely characterizes} \alephalpha, or that \alephalpha\; is
\emph{completely characterizable} by $\M$. If $<$ is a linear ordering on $\M$ and $\M$ characterizes \alephalpha, we say that $\alephalpha$ is characterizable by the  linear ordering $(\M,<)$. If $<$ is a linear ordering on a subset of $\M$, then we will say that $(\M,<)$ characterizes $\alephalpha$, if $\phi$ characterizes $\alephalpha$ and $\phi$ has a model $\N$ where $<^\N$ is of size $\alephalpha$.

We denote by \ch, the set of all
completely characterizable cardinals.
\end{df}
For now on, we consider only completely characterizable
cardinals, and we will refer to them as just characterizable cardinals.

\begin{df}\label{defhch} If $P$ is a unary predicate symbol, we say that it is
\emph{completely homogeneous} for the \lang{}- structure \A, if
$P^{\A}=\{a|\A\models P(a)\}$ is infinite and every permutation of
it extends to an automorphism of \A.

If $\kappa$ is a cardinal, we will say that $\kappa$ is
\emph{homogeneously characterizable} by
$(\phi_{\kappa},P_{\kappa})$, if $\phi_{\kappa}$ is a complete
\lomegaone- sentence, $P_{\kappa}$ a unary predicate in the
language of $\phi_{\kappa}$ such that
\begin{itemize}
    \item $\phi_{\kappa}$ doesn't have models of power $>\kappa$,
    \item if \M\; is the (unique) countable model of
    $\phi_{\kappa}$, then $P_{\kappa}$ is infinite and completely
    homogeneous for \M\; and
    \item there is a model \A\; of $\phi_{\kappa}$ such that
    $P_{\kappa}^{\A}$ has cardinality $\kappa$.
\end{itemize}

 If $(\phi_{\kappa},P_{\kappa})$ characterize $\kappa$ homogeneously and $\M,P$ are as above, we write
 $(\M,P(\M))\models (\phi_{\kappa},P_{\kappa})$. Denote the  set of all homogeneously
characterizable cardinals by \homch. Obviously,
$\homch\subset\ch$, but the inverse inclusion fails since $\aleph_0\in\ch\setminus\homch$ (cf. \cite{Knight}).
\end{df}

\begin{theorem}\label{homchimpliesloch} If $\kappa\in\homch$, then $\kappa$ is characterized by a dense linear ordering.
\begin{proof} Let $(\phi,P)$ witness the fact that $\kappa\in\homch$ and $\M$ be a countable model such that $(\M,P(\M))\models (\phi,P)$. Extend $\lang{}(\phi)$, the language of $\phi$, to include a new binary symbol $<$ and consider the new sentence $\phi'$ which is the conjunction of $\phi$ together with the sentence \[<\mbox{  is a dense linear order  on $P(\M)$ without endpoints.}\]

Now, let $\M_1,\M_2$ be two countable models of $\phi'$. Since the reducts of $\M_1,\M_2$ on the language of $\phi$ both satisfy $\phi$, they must be isomorphic. Call $i$ such an isomorphism between $\M_1$ and $\M_2$ and let $f$ be any bijection that maps the elements of  $P(\M_1)$ to the elements of $P(\M_2)$. Then $i^{-1} \circ f$ is a permutaion of $P(\M_1)$ and by the homogeneity of $P$, it extends to an automorphism of $\M_1$, say $j$. Then $i\circ j$ is a  $\lang{}(\phi)$- isomorphism between $\M_1$ and $\M_2$  and $i\circ j$ agrees with $f$ on $P(\M_1)$. In other words, there is an  $\lang{}(\phi)$- isomorphism between $\M_1$ and $\M_2$  that extends any bijection $f$ between  $P(\M_1)$ and  $P(\M_2)$. By the usual back-and-forth argument there is a bijection $f$ between $(P(\M_1),<)$ and $(P(\M_2),<)$ that preserves $<$ and by extending this $f$ we get an  isomorphism betwenn $\M_1$ and $\M_2$ that preserves both  $\lang{}(\phi)$ and $<$.
\end{proof}
\end{theorem}

Notes:(1) The assumption $\kappa\in\homch$ is too strong, since by theorem \ref{Landraitistheorem} and corollary \ref{landraitiscorollary} $\alephs{0}$ is characterizable by a linear order, while $\alephs{0}\not\in\homch$. (2) The models of $\phi'$ embed any dense linear ordering without endpoints of size up to $\kappa$.

In \cite{Knight} Hjorth proves that every $\kappa^+$ is characterized by a dense linear ordering, for all $\kappa\in\ch$\footnote{Hjorth is actually proving the result for $\kappa=\alephs{\alpha}$ and $\alpha$ countable, but his proof generalizes (cf. \cite{Knight}, proof of theorem 5.1).}.

\begin{theorem}(Hjorth) \label{loofsizekappa+} If $\kappa\in\ch$, then at least one of the following is the case:
\begin{enumerate}
  \item $\kappa^+\in\homch$ or,
  \item there is a countable model $\M$ in a language that contains a unary predicate $P$ and a binary predicate $<$ whose Scott sentence $\phiM$
  \begin{enumerate}
    \item has no models of cardinality $\kappa^{++}$ and
    \item  $\phiM$ has a model $\N$ where $(P^\N,<^\N)$ is a dense linear ordering without endpoints, it has size $\kappa^+$, and every initial segment of this linear ordering has size $\kappa$.
  \end{enumerate}
  \end{enumerate}
\end{theorem}

Cases 1 and 2 need not be exclusive the one to the other, but in either case we get

\begin{corollary}\label{loforkappaplus} If $\kappa\in\ch$, then there is a countable dense linear ordering whose Scott sentence
\begin{enumerate}
\item does not have any models in cardinality $\kappa^{++}$, but
\item does have a model which is a dense linear ordering with an increasing sequence of size $\kappa^+$.
\end{enumerate}
\begin{proof} By theorems \ref{loofsizekappa+} and \ref{homchimpliesloch}.
\end{proof}

\end{corollary}
Using theorem \ref{loofsizekappa+} we can also conclude that
\begin{corollary}\label{countablebetalo} If $\alephs{\alpha}\in\ch$, then $\alephs{\alpha+\beta}$ is characterized by a dense linear ordering, for all $\beta<\omegaone$.

In particular, $\alephs{\beta}\in\ch$ is characterized by a dense linear ordering, for all countable ordinals $\beta$.
\end{corollary}

The importance of characterizing cardinals by linear ordering is emphasized by the next theorem. It is theorem 35 from \cite{characterizablecardinals}
\begin{theorem}\label{loimplies2k} Let $\phi$ be a complete sentence such that
\begin{enumerate}
  \item For every model $\M$ of $\phi$, $<^{\M}$ is a linear order.
  \item $\phi$ does not have any models of cardinality $\lambda^+$.
  \item $\phi$ has a model $\M$ with an  $<^{\M}$- increasing sequence of size $\lambda$.
\end{enumerate}
Then $2^\lambda$ is characterizable.
\end{theorem}

Next we describe briefly a Fraisse-type construction which we are going to use.
\begin{df}Let $\A$ be a structure that contains $\M$ and if $A_0\subset \A$, then let $<A_0>$ be the substructure of $\A$ that is generated by $A_0$. We call finitely generated over $\M$ the substructures of $\A$ that have the form $<A_0>\cup\M$, where $A_0$ is a finite subset of $\A\setminus\M$. We write finitely generated/$\M$.

If $B_0=<A_0>\cup\M,\; B_1=<A_1>\cup\M$ are finitely generated/$\M$ substructures of $\A$, we write $B_0\subset B_1$ and we say that $B_0$ is a \emph{substructure} of $B_1$, if the same is true (in the usual sense) for $<A_0>$ and $<A_1>$.
\end{df}
It is straightforward to extend the above definition in the case
were we have finitely many $\M_0,\ldots,\M_n$.

Fraisse's theorem hold even for ``finitely generated/$\M$" substructures (For a proof of Fraisse's theorem one can consult \cite{shortermodeltheory}).
\begin{theorem}\label{Fraisse}(Fraisse) Fix a countable model $\M$ and let $K(\M)$ be a countable collection of finitely generated/$\M$ substructures (up to isomorphism). If $K(\M)$ has the Hereditary Property (HP), the Joint Embedding Property (JEP) and the  Amalgamation Property (AP), then there is a countable structure $\F$ which we will call the \emph{Fraisse limit} of $K(\M)$, such that
\begin{enumerate}
\item $\F$ is unique up to isomorphism and contains $\M$,
\item $K(\M)$ is the collection of all finitely generated/$\M$ substructures of $\F$ (up to isomorphism), and
\item every isomorphism between finitely generated/$\M$ substructures of $\F$ extends to an automorphism of $\F$.
\end{enumerate}

The converse is also true, i.e. if $\F$ is a countable structure such that  every isomorphism between finitely generated/$\M$ substructures of $\F$ extends to an automorphism of $\F$, and $K(\M)$ is the collection of all finitely generated/$\M$ substructures of $\F$, then $K(\M)$ has the HP, the JEP and the AP.
\end{theorem}

\begin{theorem}\label{Fraisse2}(Fraisse) Fix a model $\M$. Assume that $\A,\B$ are two structures (not necessarily countable) that contain $\M$ and  such that
\begin{itemize}
\item for every finitely generated/$\M$ substructures $C\subset D$ of $\A$ (or of $\B$), and every embedding $f:C\mapsto \A$ ($f:C\mapsto \B$), there is an embedding $g:D\mapsto \A$ ($g:D\mapsto \B$) that extends $f$, and
\item the collection of  all finitely generated/$\M$ substructures of $\A$ is the same as the collection of  all finitely generated/$\M$ substructures of $\B$.
\end{itemize}
Then $\A$ and $\B$ are back-and-forth equivalent, or $A\equiv_{\infty,\omega}\B$.

\end{theorem}

In the cases we will work with JEP follows from AP. We now give a slightly different version of the above theorem that will be more fitting to work with
\begin{theorem} \label{propertiesIandII} Fix a countable model $\M$. If $K(\M)$ is a countable collection of finitely generated/\M\; substructures (up to isomorphism) and $K(\M)$ has the HP, the JEP and the AP, then there is a unique (up to isomorphism) countable structure $\F$ that contains $\M$ and satisfies the conjunction of
\begin{description}
  \item[(I)] Every finitely generated/$\M$ substructure of $\F$ is in $K(\M)$.
  \item[(II)] For every  $A_0$ finitely generated/$\M$ substructure of \F, if there exists some $A_1\in K(\M)$ such that $A_1\supset\A_0$, then there exists some finitely generated/$\M$ substructure $B\subset\F$ and an isomorphism $i:B\cong A_1$, such that $A_0\subset B$ and $i|_{\A_0}=id$.
\end{description}
Moreover, if there is some $\lomegaone$ sentence $\psi$ such that $A\in K(\M)$ iff $A\models\psi$ (as it will the case in our example), then the conjunction of (I) and (II) can be written as a \lomegaone-sentence which is equivalent to the the Scott sentence of $\F$ and hence, it is complete.
\begin{proof} Existence follows from theorem \ref{Fraisse}. If $\F_1$ and $\F_2$ are both countable structures that contain $\M$ and satisfy $(I)$ and $(II)$, then a standard back-and-forth argument establishes the isomorphism of $\F_1$ and $\F_2$.
\end{proof}
\end{theorem}

\textbf{Notation: }In case we want to indicate which class we are talking about, we will write $(I)_{K(\M)}$ and $(II)_{K(\M)}$.

\begin{corollary}If $\M$ is countable and  $\M'\cong \M$, then $\lim K(\M')\cong \lim K(\M)$.
\end{corollary}

Theorem \ref{propertiesIandII} can be extended even in the case which $\M$ and $K(\M)$ have cardinality $\kappa>\alephs{0}$. The existence of $\F$  in this case follows from the same diagonal argument as in the countable case, but the uniqueness of the Fraisse limit fails. However,  all models of $(I)_{K(\M)}$ and $(II)_{K(\M)}$ will be $\equiv_{\omegaone,\omega}$-equivalent to each other (by theorem \ref{Fraisse2}). So, we get the following
\begin{theorem} Let $\psi$ be an $\lomegaone$ sentece. Assume that $\M$ is a countable model with Scott sentence $\phi$ and $\N$ is a model of $\phi$ (possibly uncountable) and let $K(\M)$ be the collection of all finitely generated/$\M$ substructures that satisfy $\psi$ and  let $K(\N)$ be the collection of all finitely generated/$\N$ substructures that also satisfy $\psi$. Moreover, assume that $K(\M)$ and $K(\N)$ both have the HP, the JEP and the AP. Then any model of $(I)_{K(\N)}$ and $(II)_{K(\N)}$ is $\equiv_{\infty,\omega}$- equivalent to $\lim K(\M)$.
\begin{proof} By theorem \ref{propertiesIandII}, $\lim K(\M)$ exists and it is unique, and by the comments above, $K(\N)$ has a limit which satisfies  $(I)_{K(\N)}$ and $(II)_{K(\N)}$, but this limit may not be unique.

Since $\M$ and $\N$ satisfy the same Scott sentence $\phi$, they are back-and-forth (or $\equiv_{\infty,\omega}$-) equivalent. For this it follows that for any substructure $\A\in K(\M)$ there is an  substructure $\B\in K(\N)$ such that $\A$ and $\B$ are back-and-forth equivalent, and vice versa. Using $(I)$ and $(II)$, for both $K(\M)$ and $K(\N)$, we can establish a back-and-forth equivalence for the Fraisse limits and this finishes the proof.
\end{proof}
\end{theorem}

This proves that any Fraisse limit of $K(\N)$ satisfies the Scott sentence of $\lim K(\M)$. If $\M$ is a countable model whose Scott sentence $\phi$ characterizes a certain cardinal $\kappa$, we will use the Scott sentence of $\lim K(\M)$ to characterize some cardinal $\lambda\ge\kappa$. In order to construct a Fraisse limit of $K(\N)$ we will use

\begin{theorem} \label{towardsfullness} Assume that $\M$ is a countable model whose Scott sentence $\phi$ characterizes an infinite cardinal $\kappa$, \N\; is a model of $\phi$ of cardinality $\le\kappa$, $K(\M)$ and $K(\N)$ are as above and $\lambda\ge\kappa$. Moreover, assume that:
\begin{enumerate}
  \item If A is a finitely generated/\N\; structure, then there are $\le\lambda$ many (non-isomorphic) structures in $K(\N)$ that extend A, and
  \item If $\G$ is a structure such that
     \[\N\subset\G,\; |\G\setminus\N|\le\lambda,\; \G\mbox{ satisfies }(I)_{K(\N)}\]

and for any $A_0,A_1$ are finitely generated/\N\; structures with \[A_0\subset\G,\; A_1\supset A_0\mbox{ and }A_1\in K(\N),\]

then there is another  structure $\G'$ that extends $\G$ and
\[|\G'\setminus\N|\le\lambda,\; \G'\mbox{ satisfies }(I)_{K(\N)}\]

and there is some finitely generated/\N\; structure $B\subset\G'$ and an isomorphism $i:B\cong A_1$, with $A_0\subset B$ and $i|_{A_0}=id$.
\end{enumerate}
Under the assumptions 1 and 2, we conclude that there is a structure $\G^*$ with $\N\subset\G^*$, $|\G^*|=\lambda$ and $\G^*$ satisfies $(I)_{K(\N)}$ and $(II)_{K(\N)}$. Then $\G^*$ also satisfies the Scott sentence of $\lim K(\M)$.
\begin{proof} We construct $\G^*$ by a diagonal argument. If $\G_\alpha$ is the structure at step $\alpha$ and $\G_\alpha\setminus\N$ has size $\le\lambda$, then by assumption 1, there are $\le\lambda$ many structures in $K(\N)$ that extend some finitely generated/$\N$ substructure  of $\G_\alpha$. Using the second assumption  we can ensure that we include a copy of each one of them into some $\G_\beta$, for $\beta>\alpha$.
\end{proof}
\end{theorem}

The following is theorem 27 from \cite{characterizablecardinals}
\begin{theorem}\label{ltoomegahom} If $\lambda\in\ch$, then $\ltoomega\in\homch$.
\end{theorem}

\begin{theorem}\label{cardinalitylekappa}If $\kappa$ is a cardinal in $\homch$ and $P$ is a unary predicate in a countable language $\lang{}$, then there is a sentence $\chi$ in $\lomegaone$ such that if $\N\models \chi$ then $\alephs{0}\le|P(\N)|\le\kappa$.

Moreover, if $\M$ is a countable $\lang{}$- model with $|P(\M)|=\alephs{0}$ and $$\lang{}\cap\lang{}(\chi)={P},$$ where $\lang{}(\chi)$ is the language of $\chi$, then there is a countable $\lang{}\cup\lang{}(\chi)$- model $\M'$ that extends $\M$ and $\M'$ satisfies $\chi$.

If $\N$ is a countable model isomorphic to $\M$ and $\N'$ and $\M'$ are the $\lang{}\cup\lang{}(\chi)$- extensions of $\N$ and $\M$ respectively that satisfy $\chi$, then $\M$ with the additional $\lang{}\cup\lang{}(\chi)$- structure that inherits from $\M'$ and $\N$ with the $\lang{}\cup\lang{}(\chi)$- structure that inherits from $\N'$ are isomorphic. I.e. the structure that is added on $\M$ by $\M'$ is unique (up to isomorphism).

We will say that $\chi$ witnesses the fact that  $|P(\cdot)|\le\kappa$ can be expressed in $\lomegaone$.
\begin{proof}
Let $\M_0$ be a countable model with a homogeneous predicate $P(\M_0)$ that witnesses that $\kappa\in\homch$. Let $\phi$ be the Scott sentence of $M_0$ and let $M(\cdot)$ be a unary predicate not in $\lang(\phi)$. Then take $\chi$ to be the conjunction of:
\begin{enumerate}
        \item $M(\cdot)$ and $P(\cdot)$ are disjoint,
	\item $P(\cdot)$ is infinite,
        \item $M(\cdot)\cup P(\cdot)\models \phi$ , and
	\item  $P(\cdot)$ is the homogeneous predicate of $M(\cdot)\cup P(\cdot)$ (cf. definition \ref{defhch}).
    \end{enumerate}

Since $M(\cdot)$ together with $P(\cdot)$ satisfy $\phi$, this restricts the size of $P(\cdot)$ to at most $\kappa$.

If $\M$ is a countable model as in the assumption, then let $\M'=\M\cup (\M_0\setminus P(\M_0))$ and require that $(\M_0\setminus P(\M_0))\cup P(\M)\models\phi$. Since $P$ is a homogeneous predicate, any permutation of $P(\M)$ extends to an automorphism of the whole $(\M_0\setminus P(\M_0))\cup P(\M)$ structure and the result follows.
\end{proof}
\end{theorem}
\textbf{Note}: The above proof relies  heavily on the homogeneity of $P$. If this assumption is taken away it is possible for two different extensions $(\M_0\setminus P(\M_0))\cup P(\M)$ not to be isomorphic.

The following theorem is from  \cite{Landraitis}. The interested reader should look there for more details.
\begin{theorem}(Landraitis)\label{Landraitistheorem} Let $(\M,<,P_i)_{i\in\omega}$ be a countable linear ordering, $P_i$ be a unary predicate for all $i$, and $\phi_\M$ be the Scott sentence of $\M$. Then:
\begin{enumerate}
  \item $\phi_\M$ does not have any uncountable models iff every orbit of $\M$ is scattered, and
  \item $\phi_\M$ has a model of any cardinality iff $\M$ has a self-additive interval, and
  \item if neither case (1) nor case (2) happens, then $\phi_\M$ has models in all cardinalities $\le 2^{\aleph_0}$, but no model in any cardinality above $2^{\aleph_0}$.
\end{enumerate}
All three cases do occur.
\end{theorem}

\begin{corollary}\label{landraitiscorollary} There is a countable linear ordering $(\M,<,P_i)_{i\in\omega}$ which characterizes $\alephs{0}$, i.e. the Scott sentence of $\M$ has no uncountable models.
\end{corollary}
With the above theorem Landraitis gives a complete characterization of all infinite cardinals characterized by linear orderings, but he works under the assumption that the language contains  only unary predicate symbols. Our results do not use this restriction.

\section{Dense linear orderings}

In this section we provide some background definitions and theorems about dense linear orderings that we will use later. Most of the material here follows \cite{BaumgartnerThesis}. The reader who is familiar with it can skip to the next section.

\begin{df} For infinite cardinals $\kappa\le\lambda$, let $D(\kappa,\lambda)$ iff there is a linear ordering of size $\lambda$ with a dense set of size  $\kappa$ and we let $D(\kappa,\lambda,\mu)$ iff there is a  there is a linear ordering of size $\lambda$, character $\mu$ (see definition \ref{characteroflo} below for that) and with a dense set of size  $\kappa$.

Let \[Ded(\kappa) = \sup \{\lambda|D(\kappa,\lambda) \mbox{ holds}\}\] and
\[Ded(\kappa,\mu) = \sup \{\lambda|D(\kappa,\lambda,\mu) \mbox{ holds}\}.\]
\end{df}

\begin{df} \label{characteroflo} For a linear order $(M,<)$ and some $m\in M$, the \emph{left character} of $m$ is the least cardinal $\kappa$ such that there is a cofinal function from $\kappa$ to $\{n|n<m\}$, and the \emph{right character} of $m$ is the least $\kappa$ such that there is a coinitial function from $\kappa$ to $\{n|n>m\}$.

The \emph{character} of $m$ denoted $\chi(m)$  is the least of the left and right character.

The \emph{character} of $(M,<)$ denoted $\chi(M,<)$ is the least of the cardinals $\{\chi(m)|m\in M\}$. If $(M,<)$ is a dense linear ordering, then $\chi(M,<)$ will always be infinite.
\end{df}
\begin{df}If $(M,<)$ is linear order and $(L,U)$ is a partition of $M$ with the property
\[\forall l\in L\forall u\in U,\;\; x<y,\]
then $(L,U)$ is called a \emph{Dedekind cut}. If neither $L$ has a supremum, nor $U$ an infimum, then the cut is also called a \emph{gap}.
\end{df}

\begin{df} A linear ordering $(M,<)$ is \emph{complete} if for every non-empty $M_0\subset M$ that has an upper bound, $M_0$ has a least upper bound.

The dense linear ordering $(\overline{M},\prec)$ is a \emph{completion} of the dense linear ordering $(M,<)$ if $\overline{M}$ contains $M$, $\prec$ is an extension of $<$ and $M$ is dense in $\overline{M}$, in the sense that for every $x,y\in\overline{M}$ with $x<y$, there is some $m\in M$ such that $x< m< y$.
\end{df}
The completion of a linear order is unique up to isomorphism and it is easy to see that $\chi(M,<)=\chi(\overline{M},\prec)$.
We can redefine $D(\kappa,\lambda)$ (and $D(\kappa,\lambda,\mu)$) using completions:
\begin{itemize}
\item $D(\kappa,\lambda)$ holds  iff there is a linear order of size $\kappa$ whose completion has size $\ge \lambda$ and
\item $D(\kappa,\lambda,\mu)$ holds  iff there is a linear order of size $\kappa$ and character $\mu$ whose completion has size $\ge \lambda$.
\end{itemize}

The following theorems are from \cite{BaumgartnerThesis}:
\begin{theorem} Let $\mu\le\kappa\le\lambda$ and $\mu$ is regular.
\begin{itemize}
\item [(a)] $D(\kappa,\lambda,\mu)$ holds iff there is a tree of height $\mu$ and cardinality $\le\kappa$ with at least $\lambda$ branches of length $\mu$.
\item[(b)] $D(\kappa,\lambda)$ holds iff there is a tree of height $\le \kappa$ and cardinality $\le \kappa$ with at least $\lambda$ branches.
\item [(c)] Assume that $\lambda=\sup_{\alpha<\kappa} \lambda_\alpha$. Then
\begin{eqnarray*}
\forall\alpha<\kappa\; D(\kappa,\lambda_\alpha,\mu)&\Rightarrow& D(\kappa,\lambda,\mu)\\
\forall\alpha<\kappa\; D(\kappa,\lambda_\alpha)&\Rightarrow& D(\kappa,\lambda)
\end{eqnarray*}
\end{itemize}
\end{theorem}

\begin{theorem} Let $\mu\le\kappa\le\lambda$ and $\kappa\le\kappa'$ and $\lambda'\le\lambda$. Then
\begin{eqnarray*}
D(\kappa,\lambda,\mu)&\Rightarrow & D(\kappa',\lambda',\mu)\mbox{ and }\\
D(\kappa,\lambda)&\Rightarrow &D(\kappa',\lambda')
\end{eqnarray*}
\end{theorem}

\begin{theorem}\label{nogreaterthankm} Let $\kappa,\lambda$ be infinite cardinals.
\begin{itemize}
\item If $\mu$ is the least cardinal such that $\kappa<\lambda^\mu$, then $D(\kappa,\lambda^\mu,cf(\mu))$ holds.
\item If $D(\kappa,\lambda,\mu)$ holds, then $\lambda\le \kappa^\mu$.
\end{itemize}
\end{theorem}

\begin{corollary}\label{muminimum} If $\mu$ is a regular cardinal and $\mu$ is the least such that $\kappa^\mu>\kappa$, then
$Ded(\kappa,\mu)=\kappa^\mu$.
\end{corollary}

To the best of our knowledge the following questions are open:
\begin{open} Assume that $\kappa$ is in $\ch$. Does it follow that $\kappa$ is also characterized by the Scott sentence of a dense linear ordering?

Assume that $\phi$ characterizes $\kappa$. Extend the language of $\phi$ to $\lang{1}$ by including $<$  and assume that the conjunction of  $\phi$  and ``$<$ is a dense linear ordering" has  more than $\alephs{0}$- many non-isomorphic countable models. Does it follow that one of them characterizes $\kappa$?
\end{open}

\begin{open}\label{openquestion2} Assume that $\phi$ is the Scott sentence of a dense linear ordering that characterizes $\kappa$. Can we find another dense linear ordering with Scott sentence $\psi$  such that
\begin{itemize}
\item $\psi$ characterizes $\kappa$ and
\item $\psi$ has a model with an increasing sequence of size $\kappa$?
\end{itemize}

If the answer is positive, what if we require that $\psi$ has a model with cofinality $\kappa$?
\end{open}
\begin{open} Similarly as above, if $\mu\le\kappa$, can we find some $\psi$  such that
\begin{itemize}
\item $\psi$ characterizes $\kappa$ and
\item $\psi$ has a model of size $\kappa$ and character $\mu$?
\end{itemize}
\end{open}

\section{Powers of the form $\kappa^{\alephs{1}}$} \label{kappatoalephaone}

 The main theorem in this section is
\begin{theorem}\label{kappatoaleph1} If $\kappa$ is a cardinal in $\ch$, then $\kappa^{\alephs{1}}\in\ch$.
\end{theorem}

We will prove a simingly weaker form of the theorem first

\begin{theorem}\label{kappatoaleph1least} If $\kappa$ is a cardinal in $\ch$ and $\alephs{1}$ is the least cardinal  such that  $\kappa<\kappa^{\alephs{1}}$, then $\kappa^{\alephs{1}}\in\ch$.
\end{theorem}

Before we depart in proving this theorem we make some comments:
\begin{itemize}
\item The assumption that $\alephs{1}$ is the least cardinal  such that  $\kappa<\kappa^{\alephs{1}}$, is not as restrictive as it seems.  If $\kappa<\kappa^{\alephs{0}}=\kappa^{\alephs{1}}$, then $\kappa^{\alephs{1}}$ is in $\ch$ by theorem \ref{ltoomegahom}. If $\kappa\le\kappa^{\alephs{0}}<\kappa^{\alephs{1}}$, then we can apply theorem \ref{kappatoaleph1least} to $\kappa^{\alephs{0}}$ and conclude again that $\kappa^{\alephs{1}}\in\ch$.
So, theorem \ref{kappatoaleph1least} suffices to prove theorem \ref{kappatoaleph1}.
\item If $\kappa\in\ch$ and $\alephs{1}$ is the least cardinal  such that  $\kappa<\kappa^{\alephs{1}}$, then $\kappa=\kappa^{\alephs{0}}$ and by theorem \ref{ltoomegahom},  $\kappa\in\homch$. We will make use of this fact in the proof. Also, by the same theorem $\kappa^{\alephs{1}}$ is also in $\homch$.
\item There is nothing special about $\alephs{1}$.  If $\lambda$ is a cardinal that is characterized by a dense linear ordering and $\lambda$ is the least such that $\kappa<\kappa^\lambda$, then we will prove in the next section that $\kappa^\lambda\in\ch$.
\end{itemize}

The first goal is to construct a linear order whose character is carefully controlled. In particular, we will require that it stays bounded by $\alephs{1}$. Notice that the size of the linear order will not be bounded at this point.

The idea behind the construction is to try to mimic the behavior of the lexicographic order defined on $\kappa^{\omega_1}$. For $x\neq y\in \kappa^{\omega_1}$ let $f$ be the function that maps $(x,y)$ to the least ordinal $\alpha\in\omegaone$ such that $x(\alpha)\neq y(\alpha)$. Then we can define the lexicographic order:
\[x\newlo y \mbox{ iff } x(\alpha)<y(\alpha), \mbox{ for $\alpha=f(x,y)$}.\]
Under this definition, for three distinct elements $x,y,z\in \kappa^{\omega_1}$ with $x\newlo y\newlo z$, we can have only three possibilities
\begin{enumerate}
\item Either $f(x,y)=f(x,z)=f(y,z)$, or
\item $f(x,y)=f(x,z)<f(y,z)$, or
\item $f(x,z)=f(y,z)<f(x,y)$.
\end{enumerate}
This property is the one that drives the whole construction.

By theorem \ref{loforkappaplus}, there is a linear order $(M,<)$ that characterizes $\alephs{1}$ and let $\phi$ its Scott sentence.
Let $\lang{}$ be the language that extendes the language of $\phi$ and contains the unary predicate symbols $V,M,N$, the binary predicate $\newlo$, the  binary function symbol $f$ and let  $K(\M)$ be the collection of all countable $\lang{}$- structures $\A$ that satisfies the conjunction of:
\begin{enumerate}
\item $V(\A)\cup M(\A)$ is a partition of the space. $V(\A)$ is finite, while $M(\A)$ is infinite and  $M(\A)=\M$.
\item $ \newlo$ is a linear order on $V(\A)$, not to be confused with $<$ the linear order on $M(\A)=\M$.
\item For every $x,y\in V(\A)$, $x\neq y$, $f(x,y)=f(y,x)\in M(\A)$.
\item If $x\newlo y\newlo z$ are three distinct elements of $V(\A)$, then one of the three is the case:
\begin{enumerate}
\item [4(a)] $f(x,y)=f(x,z)=f(y,z)$, or
\item [4(b)] $f(x,y)=f(x,z)<f(y,z)$, or
\item [4(c)] $f(x,z)=f(y,z)<f(x,y)$.
\end{enumerate}
\item For all $x\in \A$, $N(x)$ implies that $x\in V(\A)$. For some $x\in V(\A)$, we will say that it is $1$-colored if $N(x)$, and we will say that it is $0$-colored otherwise.
\end{enumerate}
Before we proceed we need some work. We start by the following observation on property $(4)$.

\begin{observation}\label{ob1}
(a) Property $(4)$ can be formulated equivalently as:

 If $x\newlo y\newlo z$, then $f(x,z)=\min\{f(x,y),f(y,z)\}$.

 In many cases we will use this equivalent formulation.

(b) If $f(x,y)=f(x,z)<f(y,z)$, then the only way to violate $(4)$ is if  either $y\newlo x\newlo z$, or $z\newlo x\newlo y$.
\end{observation}

\begin{df} If $\A\in K(\M)$, $x,y\in V(\A)$ and $m\in M(\A)$, define $x\simm y$ iff
\[x=y \mbox{ or } f(x,y)>m.\]
\end{df}

\begin{lemma} $\simm$ is an equivalence relation and if $m_1<m_2$, then $\sim_{m_2}$ is a refinement of $\sim_{m_1}$.
\begin{proof} Transitivity is the only one that needs some work. Assume $x\simm y$ and $y\simm z$. Then $f(x,y)>m$ and $f(y,z)>m$. Since the triplet $x,y,z$ satisfies property $(4)$, $f(x,z)$ can not be less than both $f(x,y)$ and $f(y,z)$ and the result follows.

If $m_1<m_2$ and $x\sim_{m_2} y$, then $x=y$ or $f(x,y)>m_2>m_1$.
\end{proof}
\end{lemma}

We will denote by $[x]_m$ the equivalence class of $x$ under $\simm$.

\begin{lemma} If $\A\in K(\M)$, $x,y,y'\in V(\A)$, $x\not\simm y$ and $y'\simm y$, then
\[x\newlo y\mbox{ iff } x\newlo y'.\]
\begin{proof} Since $\simm$ is an equivalence relation $x\not\simm y'$ and it suffices to prove that $x\newlo y$ implies $x\newlo y'$. Assume otherwise, i.e. $y'\newlo x\newlo y$. By property $(4)$, $f(y,y')$ is equal to at least one of the $f(x,y)$ and $f(y',x)$. Since $y'\simm y$, $f(y,y')>m$. Combining these two we conclude that either $f(x,y)>m$ or $f(y',x)>m$, which means that $x\simm y$ or $x\simm y'$. Contradiction.
\end{proof}
\end{lemma}

\begin{df} If $x,y\in V(\A)$, we write $[x]_m\newlo [y]_m$, if for some (all) $x'\in [x]_m$ and some $y'\in [y]_m$, $x'\newlo y'$.
\end{df}
In view of the above lemma, the definition is well-defined and it makes the set of all equivalence classes $\{[x]_m|x\in V(\A)\}$ into a linearly ordered set. The linear order of $V(\A)$ and the linear order on the set $\{[x]_m|x\in V(\A)\}$ are not the same, but the one arises naturally from the other, so we will use $\newlo$ for both of them.

Now we are ready to prove the following

\begin{lemma}\label{amalgamation} $K(\M)$ has the HP, the JEP and the AP.
\begin{proof}  HP is immediate and JEP follows from AP.

For AP, let $\A,\B,\C\in K(\M)$ and $\A\subset\B,\C$. We keep $N(\cdot)$ as it is on $\B$ and $\C$. All the work is to extend $\newlo$ and $f$ appropriately, so that $\B\cup \C$  becomes a structure in $K(\M)$. For all that follows fix some $b\in V(\B)\setminus V(\A)$. Let
\[m_0=\max\{f(a,b)|a\in V(\A)\},\]
\[A_0=\{a\in V(\A)|f(a,b)=m_0\},\]
\[L=\{a\in A_0|a\newlo b\} \mbox{  and  } l=\sup L,\]
\[U=\{a\in A_0|b\newlo a\} \mbox{  and  } u=\inf U.\]
Notice that $L$ or $U$ maybe empty, but at least one of them is not empty. If $L$ for instance, is empty, then let $l=-\infty$ and if $U$ is empty, then $u=\infty$. If both $L,U$ are not empty, then $f(l,u)=m_0$.

\begin{enumerate}
\item[Case 1.] If $l\neq -\infty$ and $c\in V(\C)$ with $[c]_{m_0}\newlo [l]_{m_0}$ or $[c]_{m_0}=[l]_{m_0}$, then define \[f(b,c)=\min\{f(l,c),f(l,b)\}\] and \[c\newlo b.\] If $l\newlo c$ and $f(l,c)<m_0$, then let \[f(b,c)=f(l,c)\] and \[b\newlo c.\]

\item[Case 2.] If $u\neq\infty$ we work similarly. In particular, if $[u]_{m_0}\newlo [c]_{m_0}$ or $[c]_{m_0}=[u]_{m_0}$, then define \[f(b,c)=\min\{f(u,c),f(u,b)\}\] and \[b\newlo c.\] If $c\newlo u$ and $f(u,c)<m_0$, then let \[f(b,c)=f(u,c)\] and \[c\newlo b.\]

Note: If both $l\neq -\infty$ and $u\neq\infty$, then there is some overlap between case 1 and case 2. We will prove that the definitions agree in this case (see claim \ref{wddefinition}).

\item[Case 3.] The only elements of $V(\C)$ that were not considered in the above two cases are the $c\in V(\C)$ such that $f(c,l)=m_0$ and/or $f(c,u)=m_0$, and $l\newlo c\newlo u$. Let $C_0$ be the set of all these $c$'s. If $C_0=\emptyset$, we are done, otherwise we have to do some more work:

    Split $C_0$ (arbitrarily) into two disjoint sets $D,E$ such that $D\cup E=C_0$ and for all $d\in D$ and for all $e\in E$, $d\newlo e$. Notice that we allow the possibility that one of the $D,E$ is empty. Define $d\newlo b\newlo e$, for all $d\in D$ and for all $e\in E$.
     Let $d_0=\sup D$ and $e_0=\inf E$. If both $d_0,e_0$ exist, choose one of them arbitrarily, say $d_0$, and choose some $m_1\ge f(d_0,e_0)$ and let $f(b,d_0)=m_1$ and for all other $c\in C_0$, let $f(b,c)=\min\{f(b,d_0),f(c,d_0)\}$. If only one of $d_0,e_0$ exist, say $d_0$, then choose some arbitrary $m_1\ge m_0$, let $f(b,d_0)=m_1$ and for all other $c\in C_0$, again let $f(b,c)=\min\{f(b,d_0),f(c,d_0)\}$.
\end{enumerate}

First we verify that the above definition is well-defined.

\begin{claim}\label{wddefinition} If $l\neq -\infty$ and $u\neq \infty$, then cases (1) and (2) of the above definition do not contradict each other.
\begin{proof} Notice that cases 1 and 2 overlap for all the $c$'s such that $f(l,c)<m_0$ or $f(u,c)<m_0$. If $f(l,c)<m_0$, then $f(u,c)=f(l,c)<m_0=f(l,u)$. By observation \ref{ob1}, either $c\newlo l\newlo u$, or $l\newlo u\newlo c$. Consider the first case and the second is dealt with symmetrically. If $c\newlo l$, then $[c]_{m_0}\newlo [l]_{m_0}$ and case 1 of the definition gives $f(b,c)=\min\{f(l,c),f(l,b)\}=\min\{f(l,c),m_0\}\}=f(l,c)$ and $c\newlo b$. For the same $c$, case 2 of the definition gives $f(b,c)=f(u,c)$ and $c\newlo b$. Since $f(u,c)=f(l,c)$, the two definitions completely agree.
\end{proof}
\end{claim}

Next we have to verify that $\B\cup\C$ under the above definition satisfies property $(4)$. The proof splits into many cases given by corresponding claims. We deal only with the case that $l\neq -\infty$. We can prove similar claims for the case that $u\neq \infty$, but are quite similar and we leave the details to the reader. So, for all the following claims assume that $l\neq -\infty$.

\begin{claim} \label{prelclaim2} If $c_0,c_1\in C_0$, then $f(c_0,c_1)\ge m_0$.
\begin{proof} Assume  $c_0\newlo c_1$. Then $l\newlo c_0\newlo c_1$ and $f(l,c_0)=f(l,c_1)=m_0$. So, by property $(4)$, $m_0=f(l,c_0)=f(l,c_1)\le f(c_0,c_1)$.
\end{proof}
\end{claim}

\begin{claim}\label{prelclaim} If $c\in V(\C)$ are such that $f(l,c)<m_0$, then $f(b,c)=f(l,c)$ and $c\newlo l$ iff $c\newlo b$.
\begin{proof} If $l\newlo c$, then the result is immediate by case 1. If $c\newlo l$, then $[c]_{m_0}\newlo [l]_{m_0}$ and case 1 again gives $c\newlo b$ and $f(b,c)=\min\{f(l,c),f(l,b)\}=\min\{f(l,c),m_0\}=f(l,c)$, which concludes the proof.
\end{proof}
\end{claim}

\begin{claim} If $c_0,c_1\in V(\C)$ is such that $f(l,c_0),f(l,c_1)<m_0$, then the triplet $b,c_0,c_1$ satisfies property $(4)$.
\begin{proof} By claim \ref{prelclaim}, $f(b,c_i)=f(l,c_i)$ and $c_i\newlo l$ iff $c_i\newlo b$, for $i=0,1$. Then property $(4)$ for $b,c_0,c_1$ follows from the corresponding property for $l,c_0,c_1$.
\end{proof}
\end{claim}

\begin{claim} If $c_0,c_1\in V(\C)$ are such that $f(l,c_0)<m_0<f(l,c_1)$, then the triplet $b,c_0,c_1$ satisfies property $(4)$.
\begin{proof} By claim \ref{prelclaim}, $f(b,c_0)=f(l,c_0)<m_0$ and either $c_0\newlo l\newlo b$ or $l\newlo b\newlo c_0$. Since $m_0<f(l,c_1)$, it is $[c_1]_{m_0}=[l]_{m_0}$ and by case 1, $c_1\newlo b$ and $f(b,c_1)=\min\{f(l,c_1),f(l,b)\}=m_0$. If $c_0\newlo l\newlo b$, then $c_0\newlo c_1\newlo b$. Otherwise, it would be $c_1\newlo c_0\newlo l$ and by observation \ref{ob1}, $m_0<f(l,c_1)=\min\{f(l,c_0),f(c_0,c_1)\}\le f(l,c_0)<m_0$. Contradiction. Thus, in either case $c_0$ is the minimum or the maximum of the three elements $b,c_0,c_1$. It suffices to prove $f(l,c_0)=f(c_0,c_1)$, because then $f(b,c_0)=f(l,c_0)=f(c_0,c_1)<m_0=f(b,c_1)$ and we have property $(4)$.

If $c_0\newlo c_1\newlo l\newlo b$ or $l\newlo c_1\newlo b\newlo c_0$, then by observation \ref{ob1}, $f(l,c_0)=\min\{f(l,c_1),f(c_0,c_1)\}=f(c_0,c_1)$, since $f(l,c_1)>m_0>f(l,c_0)$. If $c_0\newlo l\newlo c_1\newlo b$ or $c_1\newlo l\newlo b\newlo c_0$, then, by observation \ref{ob1} again, $f(c_0,c_1)=\min\{f(l,c_0),f(l,c_1)\}=f(l,c_0)$ and we are done.
\end{proof}
\end{claim}

\begin{claim} If $c_0,c_1\in V(\C)$ are such that $m_0<f(l,c_0),f(l,c_1)$, then the triplet $b,c_0,c_1$ satisfies property $(4)$.
\begin{proof}Without loss of generality assume that $c_0\newlo c_1$. By assumption $[c_0]_{m_0}=[l]_{m_0}=[c_1]_{m_0}$ and by case 1, $l\newlo c_0\newlo c_1\newlo b$ and $f(b,c_0)=f(b,c_1)=\min\{f(l,b),f(l,c_0)\}=\min\{f(l,b),f(l,c_0)\}=m_0$. It suffices to prove that $m_0<f(c_0,c_1)$. By observation \ref{ob1}, $f(l,c_1)=\min\{f(l,c_0),f(c_0,c_1)\}$ and by assumption, both $f(l,c_0),f(l,c_1)$ are greater than $m_0$. So, it must also be that  $f(c_0,c_1)>m_0$.
\end{proof}
\end{claim}

\begin{claim} If $c_0,c_1\in C_0$, then the triplet $b,c_0,c_1$ satisfies property $(4)$.
\begin{proof} Without loss of generality assume that $c_0\newlo c_1$ and
\[f(b,c_0)=\min\{f(b,d_0),f(c_0,d_0)\}\] and
\[f(b,c_1)=\min\{f(b,d_0),f(c_1,d_0)\}.\]
The proof splits into 3 cases:
\begin{enumerate}
\item $c_0,c_1\in D$. Then $c_0\newlo c_1\newlo d_0\newlo b$. By observation \ref{ob1}, $f(c_0,d_0)=\min\{f(c_0,c_1),f(c_1,d_0)\}$. If $f(c_0,d_0)=f(c_1,d_0)\le f(c_0,c_1)$, then $f(b,c_0)=f(b,c_1)\le f(c_0,d_0)\le f(c_0,c_1)$ and property $(4)$ is satisfied. If $f(c_0,d_0)=f(c_0,c_1)< f(c_1,d_0)$, then either $f(b,c_0)=f(b,c_1)=f(b,d_0)<f(c_0,d_0)=f(c_0,c_1)$, or $f(b,c_0)=f(c_0,d_0)=f(c_0,c_1)<f(b,d_0)$. In the latter case, since also $f(c_0,c_1)<f(c_1,d_0)$, we conclude $f(b,c_0)=f(c_0,c_1)<\min\{f(b,d_0),f(c_1,d_0)\}=f(b,c_1)$ and thus, in both cases property $(4)$ is satisfied.
\item $c_0\in D$ and $c_1\in E$. Then $c_0\newlo d_0\newlo b\newlo e_0\newlo c_1$ and by observation \ref{ob1}, $f(c_1,d_0)=\min\{f(d_0,e_0),f(c_1,e_0)\}\le f(d_0,e_0)\le m_1=f(b,d_0)$. By definition  $f(b,c_1)=\min\{f(b,d_0),f(c_1,d_0)\}=f(c_1,d_0)$. By observation \ref{ob1} again, $f(c_0,c_1)=\min\{f(c_0,d_0),f(c_1,d_0)\}$. If $f(c_0,c_1)=f(c_1,d_0)<f(c_0,d_0)$, then $f(c_0,c_1)=f(c_1,d_0)=f(b,c_1)\le f(b,d_0)$, which implies $f(c_0,c_1)=f(b,c_1)<\min\{f(c_0,d_0),f(b,d_0)\}=f(b,c_0)$ by definition and gives property $(4)$. If $f(c_0,c_1)=f(c_0,d_0)\le f(c_1,d_0)$, then $f(c_0,d_0)\le f(c_1,d_0)\le f(b,d_0)$. By definition $f(b,c_0)=\min\{f(b,d_0),f(c_0,d_0)\} = f(c_0,d_0) = f(c_0,c_1)\le f(c_1,d_0)$ and again property $(4)$ is satisfied.
\item $c_0,c_1\in E$. Then $d_0\newlo b\newlo c_0\newlo c_1$. As in the previous case, we can prove that $f(b,c_i)=f(c_i,d_0)$, $i=0,1$. Since $d_0\newlo c_0\newlo c_1$, by observation \ref{ob1}, $f(c_1,d_0)=\min\{f(c_0,c_1),f(c_0,d_0)\}$. So, either $f(b,c_1)=f(c_1,d_0)=f(c_0,c_1)\le f(c_0,d_0)=f(b,c_0)$, or $f(b,c_1)=f(c_1,d_0)=f(c_0,d_0)=f(b,c_0)<f(c_0,c_1)$, and in both cases property $(4)$ is satisfied.
\end{enumerate}
\end{proof}
\end{claim}

\begin{claim} If $c_0,c_1\in V(\C)$ and $f(l,c_0)=f(l,c_1)=m_0$, then the triplet $b,c_0,c_1$ satisfies property $(4)$.
\begin{proof} Without loss of generality assume that $c_0\newlo c_1$.  If $c_0,c_1\in C_0$, then the result is from the previous claim. Otherwise, we have to consider two cases:
\begin{enumerate}
\item Assume that $c_0\newlo l\newlo c_1$, i.e. $c_1\in C_0$, while $c_0\not\in C_0$. Then $c_0\newlo l\newlo b$ and by observation \ref{ob1}, $f(b,c_0)=\min\{f(l,c_0),f(b,l)\}=m_0$. By definition, $f(b,c_1)=\min\{f(b,d_0),f(c_1,d_0)\}$ and by definition again $f(b,d_0)=m_1\ge f(d_0,e_0)\ge m_0$, while by claim \ref{prelclaim2}, $f(c_1,d_0)\ge m_0$. So, $f(b,c_1)\ge m_0$. By observation \ref{ob1} for $c_0\newlo l\newlo c_1$, $f(c_0,c_1)=\min\{f(l,c_0),f(l,c_1)\}=m_0$. Overall, $f(b,c_0)=f(c_0,c_1)=m_0\le f(b,c_1)$ and we have property $(4)$.
\item Assume that $c_0\newlo c_1\newlo l\newlo b$, i.e. both $c_0,c_1\not\in C_0$. Then $[c_i]_{m_0}\newlo [l]_{m_0}$ and by definition $f(b,c_i)=\min\{f(l,c_i), f(l,b)\}=m_0$, for both $i=0,1$. By observatioin \ref{ob1}, $m_0=f(c_0,l)=\min\{f(c_0,c_1),f(c_1,l)\}\le f(c_0,c_1)$. Combining all these, $f(b,c_0)=f(b,c_1)=m_0\le f(c_0,c_1)$ which gives property $(4)$.
\end{enumerate}
\end{proof}
\end{claim}

\begin{claim} If $c_0,c_1\in V(\C)$ and $f(l,c_0)<m_0=f(l,c_1)$, then the triplet $b,c_0,c_1$ satisfies property $(4)$.
\begin{proof}  Without loss of generality assume that $c_0\newlo l$ (the other case is handled similarly). We split into two cases: \begin{enumerate}
\item $c_1\not \in C_0$. Then $c_0\newlo c_1\newlo l$. Otherwise, it would be $c_1\newlo c_0\newlo l$ and by observation \ref{ob1}, $m_0=f(c_1,l)=\min\{f(c_0,l), f(c_1,c_0)\}\le f(l,c_0)<m_0$. Contradiction.

So, $c_0\newlo c_1\newlo l$. By definition, $f(b,c_0)=\min\{f(l,b),f(l,c_0)\}=\min\{m_0,f(l,c_0)\}=f(l,c_0)< m_0$.  A similar argument proves that $f(b,c_1)=m_0$ and observation \ref{ob1} for $c_0\newlo c_1\newlo l$ implies $f(l,c_0)=\min\{f(l,c_1),f(c_0,c_1)\}$, while $f(l,c_0)<m_0=f(l,c_1)$. So, it must be $f(l,c_0)=f(c_0,c_1)<f(l,c_1)$. Combining all these, $f(b,c_0)=f(l,c_0)=f(c_0,c_1)< m_0=f(b,c_1)$ and we are done.

\item $c_1\in C_0$. Then $c_0\newlo l\newlo c_1$  and $f(b,c_0)=\min\{f(b,l), f(l,c_0)\}=f(l,c_0)$, while  $f(b,c_1)=\min\{f(b,d_0),f(c_1,d_0)\}$. By claim \ref{prelclaim2}, $f(c_1,d_0)$ must be greater or equal to $ m_0$ and by definition again $f(b,d_0)\ge f(d_0,e_0)\ge m_0$, which combined gives $f(b,c_1)\ge m_0$. Observation \ref{ob1} for $c_0\newlo l\newlo c_1$ gives $f(c_0,c_1)=\min\{f(l,c_0),f(l,c_1)\}=f(l,c_0)=f(b,c_0)<m_0\le f(b,c_1)$ and this concludes the claim.
\end{enumerate}
\end{proof}
\end{claim}

\begin{claim} If $c_0,c_1\in V(\C)$ and $f(l,c_0)=m_0<f(l,c_1)$, then the triplet $b,c_0,c_1$ satisfies property $(4)$.
\begin{proof} Since $f(l,c_1)>m_0$, $c_1\sim_{m_0} l$ and $[c_1]_{m_0}=[l]_{m_0}$ and by definition $c_1\newlo b$ and $f(b,c_1)=\min\{f(l,c_1),f(l,b)\}=\min\{f(l,c_1),m_0\}=m_0$. By property $(4)$ for $c_0,l,c_1$, $f(c_0,c_1)=f(l,c_0)=m_0< f(l,c_1)$.

If $c_0\newlo l\newlo b$, then by definition, $f(b,c_0)=\min\{f(l,c_0),f(l,b)\}=m_0$ and $f(b,c_0)=f(c_0,c_1)=f(b,c_1)=m_0$ and we have the result.
If $l\newlo c_0$, then either $c_0\in C_0$, or otherwise $u\neq \infty $ and $u\newlo c_0$. If $c_0\in C_0$, then $f(b,c_0)=\min\{f(b,d_0),f(c_0,d_0)\}$ and both $f(b,d_0),f(c_0,d_0)$ are $\ge m_0$. Therefore, $f(b,c_0)\ge m_0$ and $f(c_0,c_1)=m_0=f(b,c_1)\le f(b,c_0)$,which gives property $(4)$. If $u\neq \infty $ and $u\newlo c_0$, then by definition $f(b,c_0)=\min\{f(u,c_0),f(b,u)\}=m_0$ and then, $f(b,c_0)=m_0=f(b,c_1)=f(c_0,c_1)$, which concludes the proof of the claim.
\end{proof}
\end{claim}

The above claims prove that  $\B\cup\C$ satisfies property $(4)$ in all cases and the proof of the lemma is also concluded.
\end{proof}
\end{lemma}

\textbf{Note:} The proof of the above lemma would have been simpler, if we had defined $f$ and $\newlo$ differently on $C_0$. An easy example is to let $b\newlo c$, for all $c\in C_0$. The reason we went through all this work is because in the proof of theorem \ref{removinggaps} we will need a similar construction, namely we will need to make a Dedekind cut on $C_0$ and place $b$ appropriately in the cut. The details of the proof of theorem \ref{removinggaps} follow closely the proof we just did and we will omit it as it tends to be very repetitive. The interested reader should be able to fill in all the details following the example of the proof above.

Now, by theorem \ref{propertiesIandII} there is a Fraisse limit of $K(\M)$, call it $\F$ and let $\phi_\F$ be its Scott sentence. By the same theorem we know that the $\phi_\F$ is equivalent to the conjunction of $(I)_{K(\M)}$ and $(II)_{K(\M)}$ and the following lemma gives us another equivalence.

\begin{lemma}\label{Iand1and2} If $\phi_\F$ is the Scott sentence of the Fraisse limit of $K(\M)$, then $\phi_\F$ is equivalent to the conjunction of the following:
\begin{enumerate}
\item [] $(I)_{K(\M)}$ (cf. theorem \ref{propertiesIandII})
\item [$(1^*)$] for all $x\in V(\F)$ and for all $m\in M(\F)$, there exist $z_1,z_2$ such that $z_1\newlo x\newlo z_2$,  $f(x,z_1)=f(x,z_2)=m$ and we can require $z_1,z_2$ to be $0$-, or $1$- colored (not necessarily the same). In particular, $\newlo$ restricted to the subset that contains the $0$- colored (respectively the $1$- colored) elements of $V(\F)$ is a linear ordering without endpoints.
\item [$(2^*)$] for all $x\newlo y\in V(\F)$ and for all $m\ge f(x,y)$, $m\in M(\F)$, there exist some $z_1,z_2$ such that $x\newlo z_1,z_2 \newlo y$, $f(x,z_1)=m$ and $f(y,z_2)=m$ and again, we can require $z_1,z_2$ to be $0$-, or $1$- colored.
\end{enumerate}
Also note that $z_1$ is not required to be different than $z_2$, although in most cases they will be different.
\begin{proof} By theorem \ref{propertiesIandII}, $\phi_\F$ is equivalent to the conjunction of $(I)_{K(\M)}$ and $(II)_{K(\M)}$. So we have to prove that $(I)_{K(\M)}$ and $(II)_{K(\M)}$ are equivalent to the conjunction of $(I)_{K(\M)}$ and $(1^*)$ and $(2^*)$.

The direction $(I)_{K(\M)}\wedge (II)_{K(\M)}\rightarrow (I)_{K(\M)}\wedge (1^*)\wedge (2^*)$ is immediate, because if $x\newlo y\in V(\F)$, the structure generated by $x,y$ is in $K(\M)$ and then we can use $(II)_{K(\M)}$ to extend this structure to a structure that contains the desired $z_1,z_2$ for both $(1^*)$ and $(2^*)$.

The direction $(I)_{K(\M)}\wedge (1^*)\wedge (2^*)\rightarrow (I)_{K(\M)}\wedge (II)_{K(\M)}$  needs some more work. Assume that $\A_0$ is a finitely generated$/\M$ substructure of $\F$ and $A_1\supset\A_0$ with $\A_1\in K(\M)$. We need to find a finitely generated/$\M$ substructure $\F_1$ of $\F$ with $\A_0\subset\F_1$ and some isomorphism $i:\F_1\cong\A_1$ with $i|_{\A_0}=id$. We work by induction on $n=|\A_1\setminus \A_0|=|V(\A_1)\setminus V(\A_0)|$.

If $n=0$, the result is obvious and it suffices to prove the result for $n=1$. Let  $V(\A_1)=V(\A_0)\cup \{a\}$ and we will find some element $z\in\F$ such that $\A_1\cong\A_0\cup\{z\}$. If $a\newlo a_0$, or $a_0\newlo a$, for all $a_0\in V(\A_0)$, we can find this $z$ using $(1^*)$. Otherwise let $a_0$ be the maximum element in $V(\A_0)$ such that $a_0\newlo a$ and $a_1$ be the minimum element in $V(\A_0)$  such that $a\newlo a_1$, and we can find the desired $z$ using $(2^*)$ on the tuple $a_0,a_1$.
\end{proof}
\end{lemma}

Before we prove anything else we prove the following

\begin{theorem}\label{character} If $\G$ is a model of $\phi_\F$, then $(V(\G),\newlo)$ is a linear order with character $\chi(V(\G),\newlo)= cf(M(\G),<)$.
\begin{proof} Since $M(\G)$ is a model of $\phi$, it has size $\le\alephs{1}$ and the same is true for its cofinality.  Without loss of generality we will assume that $cf(M(\G),<)=\alephs{1}$. Fix a cofinal sequence $\{m_\alpha\in M(\G)|\alpha<\omegaone\}$ of length $\alephs{1}$.
It is not hard to see that
\begin{equation}\label{density}
\F\models \forall a\in V(\F)\;\forall m\in \M\; \exists b\in V(\F) \; \left(b\newlo a\wedge f(a,b)=m\right)
\end{equation}
Therefore, the same sentence is true for $\G$.

Fix some $a\in V(\G)$. For every $\alpha\in \alephs{1}$, find some $b_\alpha\newlo a$ given by (\ref{density}), such that $f(a,b_\alpha)=m_\alpha$. Then the sequence $(b_\alpha|\alpha<\alephs{1})$ is cofinal in $\{c|c\newlo a\}$. If $c\newlo a$ and $f(a,c)=m$, there is some $\alpha$ such that $m\le m_\alpha$. Then $f(b_{\alpha+1},c)=min\{f(a,c),f(a,b_{\alpha+1})\}=m$, which implies that $c\newlo b_{\alpha+1}$. So, the character of $V(\G)$ is at most $\alephs{1}$.

On the other hand, if $(c_n)_{n\in\omega}$ is an $\newlo$- increasing countable sequence, such that
\[c_0\newlo c_1\newlo\ldots\newlo a,\]
then by observation \ref{ob1}, $f(c_n,a)=\min\{f(c_{n+1},a), f(c_n,c_{n+1})\}\le f(c_{n+1},a)$. Therefore the sequence $(f(c_n,a)|n\in\omega)$ increasing and by the assumption on the cofinality of $(M(\G),<)$, this sequence can not be cofinal. Hence, it is bounded above and we can find as above some $\alpha\in\alephs{1}$ such that $c_n\newlo b_\alpha$, for all $n$. Then $b_\alpha$ is an upperbound of $(c_n))_{n\in\omega}$ and $(c_n)_{n\in\omega}$ can not be cofinal in $\{c\in V(\G)|c\newlo a\}$, which proves that the left character of (any) $a\in V(\G)$ is equal to $\alephs{1}$. We can repeat the proof for the right character being equal to $\alephs{1}$ and this concludes the proof of the theorem.
\end{proof}
\end{theorem}

\begin{lemma}\label{inflemma} If $\G$ is a model that satisfies $(I)_{K(\M)}$ and $x\in V(\G)$ and $Y\subset V(\G)$ are such that for all $y\in Y$, $y\newlo x$ and for all $m\in M(\G)$ there exists $y\in Y$ such that $f(x,y)> m$, then $x$ is the supremum of $Y$.

Quite symmetrically, if $x\newlo y$, for all $y\in Y$, and for all $m\in M(\G)$ there exists some $y$ such that $f(x,y)>m$, then $x$ is the infimum of $Y$.
\begin{proof} Towards contradiction, assume there is some $z$ such that for all $y\in Y$, $y\newlo z\newlo x$ and let $m=f(x,z)$. By assumption find some $y\in Y$ such that $f(x,y)>m$. Then the triplet $x,y,z$ contradicts property $(4)$ of $K(\M)$ since $y\newlo z\newlo x$ and $f(y,x)>m=f(x,z)$.

For the second part, the proof is symmetrical.
\end{proof}
\end{lemma}

\begin{df}\label{irremovable gaps} Assume $\G$ is a structure that satisfies $(I)_{K(\M)}$ and $(D,E)$ is a Dedekind cut, i.e. a partition of $V(\G)$ such that $d\newlo e$, for all $d\in D$ and for all $e\in E$. If there exists some $d\in D$ such that the set $\{f(d,e)|e\in E\}$ is cofinal in $\{f(d',e)|d'\in D,e\in E\}$ and the set $\{f(d,d')|d'\in D,d\newlo d'\}$ is coinitial in $\{m\in \M|m>f(d',e)\mbox{ for all $d'\in D,e\in E$}\}$ and the set $\{f(d',e)|d'\in D,e\in E\}$ does not have a supremum (equivalently, the set $\{m\in \M|m>f(d',e)\mbox{ for all $d'\in D,e\in E$}\}$ does not have an infimum), then we will call $(D,E)$  an \emph{irremovable} gap.

Symmetrically, if there exists some   $e\in E$ such that the set $\{f(d,e)|d\in D\}$ is cofinal in $\{f(d,e')|d\in D,e'\in E\}$ and the set $\{f(e,e')|e'\in E,e'\newlo e\}$ is coinitial in $\{m\in \M|m>f(d',e)\mbox{ for all $d'\in D,e\in E$}\}$ and  the set $\{f(d,e')|d\in D,e'\in E\}$ does not have a supremum (equivalently, the set $\{m\in \M|m>f(d',e)\mbox{ for all $d'\in D,e\in E$}\}$ does not have an infimum), then $(D,E)$ is an irremovable gap.

If $(D,E)$ consists of a gap on $V(\G)$ that is not irremovable, we will call it a \emph{removable} gap.
\end{df}

If $(D,E)$ is an irremovable gap, it follows that neither $\sup D$ exists nor $\inf E$, and moreover, we can not extend $V(\G)$ by adding an element $b$ such that for all $d\in D$ and for all $e\in E$, $d\newlo b\newlo e$. The reason for that is that if $d$ witnesses the fact that $(D,E)$ is irremovable, then by observation \ref{ob1}, $f(b,d)$ has to be smaller or equal than all the elements of the set $\{f(d,d')|d'\in D,d\newlo d'\}$ and bigger or equal to all the elements of the set $\{f(d,e)|e\in E\}$, which by tbe assumptions above can not happen.

Also notice that the sets  $$M_1=\{m\in \M|m\le f(d',e)\mbox{ for some $d'\in D,e\in E$}\}$$ and $$M_2=\{m\in \M|m>f(d',e)\mbox{ for all $d'\in D$, $e\in E$}\}$$ consist of a Dedekind cut on $(\M,<)$. If neither $M_1$ has a supremum nor $M_2$ an infimum, then the cut is a gap and this gap gives rise to irremovable gaps on $(V(\G),\newlo)$ (one for every element of $V(\G)$). Therefore, if $(\M,<)$ is a complete order,  then there are no irremovable gaps on $(V(\G),\newlo)$. Otherwise the number of irremovable gaps is bounded by $|V(\G)|$ times the number of gaps on $(\M,<)$. If the character of  $(\M,<)$ is $\mu$, $\mu\le \alephs{1}=|\M|$, then the number of gaps on $(\M,<)$ is $\le \alephs{1}^\mu\le 2^{\alephs{1}}$.

\begin{theorem}\label{removinggaps} Assume $\G$ is a structure that satisfies $(I)_{K(\M)}$, $|\G\setminus\M|\le\kappa$, and $(D,E)$ is a removable gap on $V(\G)$. Also assume that $\A_0\subset\G$, $\A_1=\A_0\cup\{b'\}$ and  $\A_1$ is a structure in $K(\M)$. Then there exists another  structure $\G'$ which extends $\G$, $\G'$ satisfies $(I)_{K(\M)}$,  $|\G' \setminus \M|\le\kappa$, there exists some element $b\in\G'$ such that $d\newlo b\newlo c$, for all $d\in D$ and all $e\in E$, and there exists an isomorphism  $i:\A_0\cup\{b\}\cong\A_1$, with $i|_{\A_0}=id$.

We will say that $\G'$ removes that gap $(D,E)$ by adding $b'$.
\begin{proof} As in lemma \ref{amalgamation}, let $d_0=\sup D$ and $e_0=\inf E$ and notice that $d_0, e_0$ may not exist. We distinguish the following cases:

\begin{itemize}
\item $d_0$ exists and there is some $m\in M(\G)$ such that for all $e\in E$, $f(d_0,e)\le m$.

Then choose some arbitrary $m_1\ge m_0$, let $f(b,d_0)=m_1$ and for all other $c\in C_0$, let $f(b,c)=\min\{f(b,d_0),f(c,d_0)\}$.

\item  Symmetrically, if $e_0$ exists and there is some $m\in M(\G)$ such that for all $d\in D$, $f(e_0,d)\le m$, then  choose some arbitrary $m_1\ge m_0$, let $f(b,e_0)=m_1$ and $f(b,c)=\min\{f(b,e_0),f(c,e_0)\}$, for all other $c\in C_0$.

\item  If both $d_0,e_0$ exist, then choose one of them arbitrarily, say $d_0$, and choose some $m_1\ge f(d_0,e_0)$ and let $f(b,d_0)=m_1$ and $f(b,c)=\min\{f(b,d_0),f(c,d_0)\}$,  for all other $c\in C_0$.

These first three cases are similar to the ones we encountered in the proof of Amalgamation. The next ones are new:

\item $d_0$ exists and for all $m\in M(\G)$ there exists $e\in E$, $f(d_0,e)\ge m$.

Then by lemma \ref{inflemma} $d_0$ is the infimum of $E$ and in this case $(D,E)$ is not a gap.

\item Symmetrically, if $e_0$ exists and for all $m\in M(\G)$ there exists $d\in D$, $f(e_0,d)\ge m$, then $e_0$ is the supremum of $D$ and again $(D,E)$ is not a gap.

\item If for every $d\in D$  the set $\{f(d,e)|e\in E\}$ is not cofinal in $\{f(d',e)|d'\in D,e\in E\}$, then for every $d$ there exists some $d'\in D$ and $e'\in E$  such that for every $e\in E$, $f(d,e)< f(d',e')$. In this case define $f(b,d)=f(d,d')$.

Similarly, if for every $e\in E$ the set  $\{f(d,e)|d\in D\}$ is not cofinal in $\{f(d,e')|d\in D,e'\in E\}$, then there exists some $e'\in E$ and some $d'\in D$ such that for every $d\in D$, $f(d,e)<f(d',e')$. Define $f(b,e)=f(e,e')$.

\item If there exists some $d\in D$ such that  the set $\{f(d,e)|e\in E\}$ is cofinal in $\{f(d',e)|d'\in D,e\in E\}$ and there exists some $s\in \M$ greater than all $f(d',e)$, $d'\in D$, $e\in E$ and smaller than all $f(d,d')$, $d'\in D$, $d\newlo d'$, then let $f(b,d)=s$ and for every other $c\in C_0$, let $f(b,c)=\min\{f(b,d),f(c,d)\}$.

\item Symmetrically, if there exists some $e\in E$ such that  the set $\{f(d,e)|d\in D\}$ is cofinal in $\{f(d,e')|d\in D,e'\in E\}$ and there exists some $s$ greater than all $f(d,e')$, $d\in D$, $e'\in E$ and smaller than all $f(e,e')$, $e'\in E$, $e'\newlo e$, then let $f(b,e)=s$ and for every other $c\in C_0$, let $f(b,c)=\min\{f(b,e),f(c,e)\}$.


\end{itemize}
In all these cases we have to prove that $f$ was defined in such a way that $\G'=\G\cup\{b\}$ also satisfies $(I)_{K(\M)}$. The details of the proof follow the proof of lemma \ref{amalgamation}  and are left to the reader.
\end{proof}
\end{theorem}

\begin{theorem}\label{embedanylo} Let $\M$ be a model of $\phi$ of cardinality $\alephs{1}$. Then for every linear ordering $(L,<)$ that has cardinality $|L|=\kappa\ge\alephs{1}$, there is some model $\G^*$  of cardinality $\kappa$ that satisfies $\phi_\F$, $M(\G^*)=\M$, and there exists a one-to-one function $F:L\rightarrow V(\G^*)$ such that
\[x< y\mbox{ iff } F(x)\newlo F(y)\] and
\[\mbox{ for all $x\in L$, $F(x)$ is $0$- colored}.\]
In other words, every linear ordering of size $\ge\alephs{1}$ can be embedded into a model of $\phi_\F$ of the same cardinality.
\begin{proof} Let $\G$ be the $\lang{}$- structure such that $(V(\G),\newlo)\cong (L,<)$ and for every $x,y\in V(\G)$, $x,y$ are both $0$- colored and $f(x,y)=m$, for some fixed $m\in M(\G)$. It is immediate that this $\G$ satisfies $(I)_{K(\M)}$ and we have to extend $\G$ to a structure $\G^*$  that satisfies both  $(I)_{K(\M)}$ and  $(II)_{K(\M)}$. Since for any finitely generated/$\M$ substructure of $\G$ there are $\alephs{1}\le\kappa$ many structures in $K(\M)$ that extend it, we can use theorem \ref{towardsfullness} and theorem \ref{removinggaps} to find the desired $\G^*$.
\end{proof}
\end{theorem}

\begin{corollary}\label{modelinanykappa} For every infinite cardinal $\kappa$, there is a  model of $\phi_\F$ that has size $\kappa$.
\end{corollary}

\begin{theorem}\label{completiontheorem} Let $\G$ be a model of $\phi_\F$ and $\overline{V(\G)}$ the $\newlo$- completion of $V(\G)$. Then there is some model $\overline{\G}$ such that
\begin{itemize}
\item If $U$ is the set of irremovable gaps of $V(\G)$ (cf. definition \ref{irremovable gaps}), then $V(\overline{\G})=\overline{V(\G)}\setminus U$,
\item for all $x\in V(\overline{\G})\setminus V(\G)$, $x$ is $1$-colored, i.e. no new $0$-colored elements are introduced,
\item the function $f^{\overline{\G}}$ restricted on $\G\times\G$ agrees with the function $f^\G$ and
\item $\overline{\G}\models \phi_\F$.
\end{itemize}
The model $\overline{\G}$ we will call the \emph{completion} of $\G$.

In particular, if  $|\G|=\kappa> 2^{\alephs{1}}$, then the completion of $\G$ has cardinality equal to $|\overline{V(\G)}|$.
\begin{proof}  Since $\G$ satisfies $\phi_\F$, by theorem \ref{Iand1and2}, $\G$ satisfies $(I)_{K(\M)}$ and $(1^*)$ and $(2^*)$. Since  $(1^*)$ and $(2^*)$ are density requirements and $V(\G)$ is dense in $V(\overline{\G})$, $\overline{\G}$ satisfies $(1^*)$ and $(2^*)$ too. So, it remains to show that $\overline{\G}$ satisfies $(I)_{K(\M)}$ and this follows from theorem \ref{removinggaps}. We remove all the removable gaps, one at a time, by applying theorem \ref{removinggaps} and we make sure that all the elements we add are $1$- colored.

If $\kappa> 2^{\alephs{1}}$, then there exist at most $\kappa\cdot \alephs{1}^\mu \le \kappa\cdot 2^{\alephs{1}}=\kappa$ many irremovable gaps (see comments after definition \ref{irremovable gaps}), where $\mu$ is the character of $(\M,<)$ and $\mu\le\alephs{1}$, and the result follows.
\end{proof}
\end{theorem}

\begin{theorem} If $\kappa\in\ch$ and $\alephs{1}$ is the least cardinal such that $\kappa^{\alephs{1}}>\kappa$, then $\kappa^{\alephs{1}}\in \ch$.
\begin{proof} If $\kappa\le 2^{\alephs{1}}$, then $\kappa^{\alephs{1}}=2^{\alephs{1}}$ and the result follows from theorem \ref{loimplies2k} and corollaries \ref{loforkappaplus} and \ref{countablebetalo}.
So assume that $\kappa> 2^{\alephs{1}}$. By the assumption on $\kappa$ and by theorem \ref{ltoomegahom} we can assume that $\kappa\in\homch$. Let $\psi$ be the conjunction of $\phi_\F$ together with the sentence that expresses the fact that the set of the $0$- colored elements of $V(\cdot)$ has cardinality $\le\kappa$. Since $\kappa\in\homch$ we can express this fact by a sentence in $\lomegaone$ by theorem \ref{cardinalitylekappa}.

By theorem \ref{embedanylo}, $\phi_\F$ has a model $\G^*$  that embeds the linear ordering $\kappa^{<\alephs{1}}$ into the set of the $0$- colored elements. By assumption $\kappa^{<\alephs{1}}=\kappa$ and we can assume that $\G^*$ also has cardinality $\kappa$. From this it follows that $\G^*$ is also a model of $\psi$. By theorem \ref{completiontheorem} there is a model $\overline{\G^*}$ of cardinality equal to the cardinality of $|\overline{\kappa^{<\alephs{1}}}|=\kappa^{\alephs{1}}$ that introduces no new $0$- colored elements. Therefore, $\overline{\G^*}$ is also a model of $\psi$. Since by theorem \ref{character} any model of $\psi$ (and $\phi_\F$) has character $\le \alephs{1}$, and any model of $\psi$ has a dense subset of cardinality $\le\kappa$ (the set of the $0$-colored elements),  by theorem \ref{nogreaterthankm} there is no model of $\psi$ of cardinality $>\kappa^{\alephs{1}}$ which concludes the proof.
\end{proof}
\end{theorem}

\begin{theorem} \label{kappatoalpheoneinch}If $\kappa\in\ch$, then $\kappa^{\alephs{1}}\in\ch$.
\begin{proof} If $\kappa^{\alephs{0}}=\kappa^{\alephs{1}}$, the result follows from \ref{ltoomegahom}. Otherwise  use the previous theorem for $\kappa^{\alephs{0}}$.
\end{proof}
\end{theorem}

\section{Powers of the form $\kappa^{\lambda}$} \label{kappamu}

There is nothing particular about $\alephs{1}$ that can not be generalized to any uncountable cardinal $\lambda$ that is characterized by a dense linear ordering.  The proofs of the following theorems follow from the proof of the corresponding theorems in the previous section by replacing $\alephs{1}$ with $\lambda$. Thus we have the following:

Let $(\M,<)$ be a countable dense linear ordering that characterizes  $\lambda$ and $K(\M)$ be defined as above (see theorem \ref{kappatoaleph1least}) . Then $K(\M)$ has the HP, the JEP and the AP and there exists a Fraisse limit for $K(\M)$ which we will call $\F$ and let $\phi_\F$ be its Scott sentence.

\begin{theorem}\label{character2} If $\G$ is a model of $\phi_\F$, then $(V(\G),\newlo)$ is a dense linear ordering with character $\chi(V(\G),\newlo)= cf(M(\G),<)$.
\end{theorem}
Since the character of $(M(\G),<)$ is bounded by $|M(\G)|\le\lambda$, the set of gaps on $(M(\G),<)$ will have size at most $2^\lambda$. If $\lambda$ is a singular cardinal, then $\chi(V(\G),\newlo)= cf(M(\G),<)<\lambda$ and we may get strict inequality.

\begin{theorem}\label{embedanylo2} If $\lambda$ and $\phi_\F$ are as above, then for every linear ordering $(L,<)$ with cardinality $|L|=\kappa\ge\lambda$, there exists some model $\G^*$ of $\phi_\F$ of cardinality $\kappa$  and a one-to-one function $F:L\rightarrow \G^*$ such that
\[x< y\mbox{ iff } F(x)\newlo F(y)\]
and
\[\mbox{ for all $x\in L$, $F(x)$ is $0$- colored}.\]

In other words, every linear ordering of size $\ge\lambda$ can be embedded into a model of $\phi_\F$ of the same cardinality.
\end{theorem}
 If $\G$ is a model of $\phi_\F$ of size $\kappa$, then the set of irremovable gaps on $V(\G)$ has size $\le\kappa\cdot 2^\lambda$.

\begin{theorem} If $\G$ is a model of $\phi_\F$ of cardinality $\kappa>2^\lambda$, then there is a model of $\phi$ of cardinality  equal to $|\overline{V(\G)}|$.
\end{theorem}

\begin{theorem}\label{leastlambda} Assume that  $\kappa>2^{\lambda}$ is a characterizable cardinal, $\lambda$ is a cardinal characterizable by a dense linear ordering and  $\lambda$ is the least cardinal such that $\kappa^{\lambda}>\kappa$. Then $\kappa^{\lambda}\in \ch$.
\end{theorem}
Observe here that we allow the possibility that $\lambda$ is countable. In this case the result follows from theorem \ref{ltoomegahom}. In the case that $\kappa\le 2^{\lambda}$, the characterizability of $\kappa^{\lambda}=2^{\lambda}$ follows from theorem \ref{loimplies2k}, but we require the extra assumption that there exists a model with an increasing sequence of size $\lambda$. If $\lambda$ is a successor of a characterizable cardinal (as it is in the case of $\alephs{1}$), such a sequence is guaranteed by corollary \ref{loforkappaplus}, but in the general case this question is open (see Open Question \ref{openquestion2}). So, we have
\begin{theorem} \label{powerclusters} If $\alephs{\alpha}$ and $\kappa^{\alephs{\alpha}}$ are both in $\ch$, then $\kappa^{\alephs{\alpha+\beta}}\in\ch$, for all countable ordinals $\beta$.
\begin{proof} By induction on $\beta$. If $\kappa^{\alephs{\alpha+\beta}}=2^{\alephs{\alpha+\beta}}$, then  use corollary \ref{loforkappaplus} and theorem \ref{loimplies2k}. Otherwise, use corollary \ref{countablebetalo} and theorem \ref{leastlambda}.
\end{proof}
\end{theorem}

\begin{corollary}\label{countable powers} If $\kappa\in\ch$, then $\kappa^{\alephs{\alpha}}\in\ch$, for all countable ordinals $\alpha$.
\end{corollary}

\begin{theorem}\label{powerclusters2} Assume that  $\kappa>2^{\alephs{\alpha}}$ is a characterizable cardinal, $\alephs{\alpha}$ is a cardinal characterizable by a dense linear ordering and  $\alephs{\alpha}$ is the least cardinal such that $\kappa^{\alephs{\alpha}}>\kappa$. Then $\kappa^{\alephs{\beta}}$ is in $\ch$, for all $\beta<\alpha+\omegaone$.
\begin{proof} By theorems \ref{powerclusters} and \ref{leastlambda}.
\end{proof}
\end{theorem}



\end{document}